\documentclass [11pt]{amsart} 
\usepackage {amsmath,amssymb, amscd, amsthm, epsfig, color}
\input xy
\xyoption {all}

\newtheorem{theorem}{Theorem} 
\newtheorem {lemma}{Lemma}
\newtheorem {lemmaa}{Lemma}[section]
 
\newtheorem {conjecture}{Conjecture}
\newtheorem {question}{Question}

\theoremstyle{definition}

\theoremstyle {definition}
\newtheorem {definition}{Definition}

\theoremstyle {definition}
\newtheorem {fact}{Fact}
\newtheorem* {claim}{Claim}

\def \pr {{\mathbb {P}^r}}
\def \f {{\mathcal F}}
\def \ms {{\mathbb S}}
\def \fstack {{\mathcal {\overline M}_{0,n}(X, \beta)}}
\def \fscheme {{\overline {M}_{0,n} (X, \beta)}}
\def \fopenscheme {{M_{0,n}(X, \beta)}}
\def \p {{\mathbb {P}}}
\def \mon {{M_{0,n}}}
\def \monbar {{\overline {M}_{0,n}}}

\def \s {{\mathcal {S}}}
\def \q {{\mathcal {Q}}}
\def \e {{\mathcal E}}
\def \o {{\mathcal O_{\mathbb P^1}}}
\def \c {{\mathbb C}}
\def \cs {{\mathbb C^{\star}}}
\def \mor {{\text {Mor}_d(\mathbb P^1, X)}}

\begin {document}
\baselineskip=16pt

\title[The Tautological Rings of the Moduli Spaces of Stable Maps]{The Tautological Rings of the Moduli Spaces of Stable Maps to Flag Varieties}
\author {Dragos Oprea}
\address {Department of Mathematics}
\address {Massachusetts Institute of Technology}
\address {77 Massachusetts Avenue, Cambridge, MA 02139.}
\email {oprea@alum.mit.edu}
\date{} 

\begin {abstract}

We show that the rational cohomology classes on the moduli spaces of genus zero stable maps to $SL$ flag varieties are tautological. 

\end {abstract}

\maketitle

The Kontsevich moduli stacks of stable {\it maps} arise as generalizations of the classical Deligne-Mumford spaces of stable {\it curves}. The intersection theory of the Kontsevich spaces has been intensively studied in the last decade in relation to enumerative geometry and string theory.

Partial results are known about the cohomology or the Chow groups of the Deligne-Mumford spaces in low codimension or low genus. Higher genera are particularly difficult. 
In this case, one may study the tautological rings, bearing in mind that nontautological classes do exist \cite {taut}. 
By contrast, the genus zero case is well understood. Keel proved that the cohomology is tautological, in fact generated by boundary classes of curves with fixed dual graph \cite {K}. A quick Hodge theoretic proof of Keel's result is presented in \cite {zero}. This result has implications, for instance, in the study of the tree-level cohomological field theories \cite {M}.

\subsection {The tautological rings.} As the moduli spaces of stable curves are examples of Kontsevich spaces, it was suggested in \cite {three} that it may be useful to push the investigation of the tautological rings in the context of Gromov-Witten theory. Here, we study the generalization of Keel's theorem to the genus zero Kontsevich spaces of 
maps to flag varieties $X$.

There are natural cohomology classes defined on the moduli spaces of stable maps; their intersection numbers can be expressed in terms of the Gromov-Witten invariants of $X$. We show that these natural classes generate the rational cohomology. This implies that the Gromov-Witten invariants essentially capture the {\it entire} intersection theory of the Kontsevich moduli spaces. 

Among the natural cohomology classes, we single out the boundary classes of maps with fixed dual graph. We may impose additional constraints making the marked points or nodes map to certain Schubert subvarieties of $X$ and requiring that the image of the map intersect various Schubert subvarieties. There is a cleaner way of bookkeeping the geometric Schubert-type classes we just described which leads to the definition of the tautological rings. 
\vskip.1in
To set the stage, we let $X=G/P$ be a homogeneous variety, and $\beta\in H_{2}(X,\mathbb Z)$ a homology class. The moduli stacks $\mathcal {\overline M}_{0,S} (X, \beta)$ parametrize $S$-pointed genus zero stable maps to $X$ in the homology class $\beta$. We use the notation $\fstack$ when the marking set $S=\{1, 2, \ldots, n\}$. Additionally, we denote by $\fscheme$ the corresponding coarse moduli scheme.

The stable map spaces are connected by a system of natural morphisms, which we enumerate below: \begin {itemize} 
\item \text {forgetful morphisms}: $\pi:\mathcal {\overline M}_{0,S} (X, \beta) \to {\mathcal {\overline M}}_{0, T}(X, \beta)$ \text {defined whenever} $T\subset S$, \item {\text evaluation morphisms to the target space}, $ev_i:\mathcal {\overline M}_{0,S} (X, \beta) \to X$ for all $i\in S,$ \item \text {gluing morphisms which produce maps with nodal domains},  $$gl:{\mathcal {\overline M}}_{0, S_1\cup \{\bullet\}}(X, \beta_1) {\times_{X}} {\mathcal {\overline M}}_{0, \{\bullet\} \cup S_2} (X, \beta_2) \rightarrow {\mathcal {\overline M}}_{0, S_1\cup S_2}(X, \beta_1+\beta_2).$$ \end {itemize}

More generally, fixing a dual graph $\Gamma$ corresponding to genus $0$ stable maps to $X$, we can consider the fibered product $\overline{\mathcal M}(\Gamma)$ of Kontsevich spaces parametrizing stable maps with the given dual graph $\Gamma$. Each factor of the fibered product $\overline{\mathcal M}(\Gamma)$ is of the form $\overline{\mathcal M}_{0, n_v}(X, \beta_v)$ for vertices $v$ of the graph $\Gamma$; see Section \ref{s1111} for more details. The stacks $\overline{\mathcal M}(\Gamma)$ are smooth Deligne-Mumford \cite{KP}. The forgetful and gluing morphisms are also defined on the spaces $\overline{\mathcal M}(\Gamma)$, for different graphs $\Gamma$. 

The classes pulled back from the target serve as the seed data for the tautological systems. We get more classes making use of the natural morphisms. The definition below imitates the more familiar one for the moduli spaces of stable curves. 

\begin {definition}
The genus $0$ tautological rings $R^{\star}(\overline{\mathcal M}(\Gamma))$  are the smallest system of {\it subrings} of the rational cohomology $H^{\star}(\overline{\mathcal M}(\Gamma))$ such that:  
\begin {itemize}
\item the system is closed under pushforwards by the natural forgetful and gluing morphisms, 
\item for all evaluation morphisms ${ev}:\overline{\mathcal M}(\Gamma)\to X$, the evaluation classes $ev^{\star} \alpha$ for $\alpha \in H^{\star} (X)$ are in the system.
\end {itemize}
\end {definition}

As a special case, we obtain the tautological rings $$R^{\star}(\fstack)\subset H^{\star}(\fstack).$$ 

Typical examples of tautological classes are the following. Let $\alpha_1, \ldots, \alpha_p$ be classes on $X$. The class $\kappa_n(\alpha_1, \ldots, \alpha_p)$ is obtained by pushforward under the forgetful projection $\pi:{\mathcal {\overline M}}_{0, n+p}(X, \beta)\to \fstack$: \begin {equation}\label {kappa}\kappa_n(\alpha_1, \ldots, \alpha_p)=\pi_{\star}(ev_{n+1}^{\star} \alpha_1 \cdot \ldots \cdot ev_{n+p}^{\star} \alpha_p ).\end {equation} The virtual fundamental cycles of various moduli spaces of stable maps to smooth subvarieties of $X$ are examples of classes which are not immediately seen to be tautological. 

The main result of this paper is the following extension of Keel's theorem from the case of rational marked curves to the case of genus $0$ stable maps to $SL$ flag varieties.  

\begin {theorem}\label {main}
Let $X$ be any $SL$ flag variety over the complex numbers. Then all rational cohomology classes on $\fstack$ are tautological. 
\end {theorem}

The result contained in the theorem above was already known for $X=\pr$ in complex codimension $1$ by work of Pandharipande \cite {divisors}. Of course, the degree $0$ case was solved by Keel. The degree $1$ case follows from work of Fulton and MacPherson \cite {FM}. Behrend and O'Halloran \cite {BH} computed the cohomology ring in degrees $2$ and $3$ for maps to $\pr$ without markings, essentially using localization techniques -- their result implies ours in degrees $2$ and $3$. In higher degrees or for general flags such a result, although not surprising, was indeed missing. After our paper was submitted, Mustata and Mustata announced a presentation of the cohomology ring for all degrees, when the target is again $\pr$ \cite {MM}.

We should mention that the topologists gained a good understanding of the spaces of holomorphic maps to flag varieties $\text {Map}_{\beta}(\mathbb P^1, X)$ by comparing them with the spaces of continuous maps; work in this direction was initiated by Segal \cite{S}. However, our results are of entirely different nature - by contrast with the spaces studied in topology, the stable map spaces have only algebraic cohomology. 

\subsection {Idea of the proof.} As it often happens, one attempts to apply localization techniques to understand the classes on the Kontsevich moduli spaces. Such an approach is tempting in our case as well, especially because the fixed point loci are known to be products of moduli spaces of rational marked curves whose cohomology groups are indeed tautological. However, the author could not obtain a full proof of the above theorem following this line of reasoning unless the target is $\pr$. We will present the argument in \cite {O1}. 

Nonetheless, localization gives the weaker result below, whose proof was communicated to us by Professor Rahul Pandharipande: \vskip.1in

{\it Let $X=G/P$ be a homogeneous space acted on by the maximal torus $T$. The equivariant tautological rings $R_T^{\star}(\fstack)$ are obtained by making use of the equivariant morphisms in Definition $1$. Then $R^{\star}_T(\fstack)$ and $H^{\star}_T(\fstack)$ become isomorphic {\it after inverting the torus characters}.}
\vskip.1in
Instead of using localization, our proof of the theorem relies on the Deligne spectral sequence as in \cite {zero}. Most commonly, this spectral sequence is used to define the mixed Hodge structure on the cohomology of smooth open varieties. We will use the reverse procedure: we will identify the lowest weight Hodge piece in the cohomology of the open stratum and use it to derive information about the cohomology of the compactified moduli spaces. The argument is inductive, making use of the fact that the open stratum is compactified by adding normal crossing divisors which are essentially lower dimensional moduli spaces of stable maps. Two observations are necessary in order to identify the relevant piece of the Hodge structure. First, the morphism spaces are hard to understand cohomologically. For this reason, we will appeal, as in \cite {BDW}, to the different compactification provided by the Hyper-Quot scheme. Secondly, we use a technique which goes back to Atiyah and Bott in their study of the moduli of bundles over curves: when the number of markings is small, we express the open stratum as a global quotient and carry out the computation in equivariant cohomology. A similar program was partially pursued in \cite {chow} for unpointed maps to $\pr$. The general case is more involved. When $X$ is a Grassmannian, we rely on Stromme's description of the $Quot$ scheme as the base of a principal bundle sitting in an affine space \cite {quot}. Combined with the Atiyah-Bott technique we obtain enough information about the Hodge structure of the open stratum to prove our main result. When $X$ is any $SL$ flag variety, the proof above does not immediately carry over. We will make use of the already proved results for Grassmannians combined with a general statement about the cohomology of the Hyper-Quot scheme $HQuot$, which we obtain using a technique due to Ellingsrud-Stromme-Beauville \cite {ES}. This will be sufficient to finish the proof. 
\vskip.1in
\noindent {\bf Conventions and remarks.}\vskip.1in \begin {itemize} \item [(i)] All stacks considered in this paper will be smooth and of Deligne Mumford type over $\c$. This allows us to speak about their Chow {\it rings} as defined in \cite {V1} and of their cohomology as defined in \cite {Be}. In this paper, {\it rational coefficients} will always be understood. 
\item [(ii)] The results of this paper hold either in cohomology or in the Chow groups. We will write most of the arguments in cohomology. However, we show in \cite {O1} that as a consequence of localization, {\it rational cohomology and rational Chow groups of the compactified spaces are isomorphic.} For the reader interested in the Chow groups of the {\it open strata}, we indicate the necessary changes in our proofs at the appropriate places. 
\item [(iii)] In cohomology, we need to keep track of the Hodge structure. To compensate, the cohomological proofs show that the classes supported on the boundary fibered products are tautological (Lemma $\ref{orbilemma}$). This is the {\it only} part of the argument which we do not carry out in the Chow groups. It is conceivable that one can write down an entirely algebraic proof of our main result\footnote {The referee remarks that the techniques of \cite {T} could possibly provide an alternative to our Hodge theoretic arguments.}, which would then work over any field, possibly exploiting the torus action. For instance, one can show in this fashion that the cohomology and Chow rings of the boundary fibered products $\overline{\mathcal M}(\Gamma)$ are isomorphic.  
\item [(iv)] We may consider the rational cohomology (or rational Chow rings) of the coarse moduli schemes. Pullback/pushforward by $p: \fstack\to \fscheme$ induce isomorphisms in cohomology (Chow) \cite {Be}, \cite {V1}. 
We identify the tautological rings of the moduli stacks and of the coarse moduli schemes via $p$. 
\item [(v)] It is customary to include the $\psi$ classes in the definition of the tautological rings of the Deligne Mumford spaces. However, in our case, Lemma $2.2.2$ in \cite {divisors} expresses the $\psi$'s in terms of evaluation classes, boundaries and the $\kappa$ classes. Here, $\psi_i=c_1(\mathcal L_i)$, where the fiber of the line bundle $\mathcal L_i \to \fstack$ over a stable map with domain $C$ and markings $x_1, \ldots x_n$ is the $i^{\text {th}}$ cotangent line $T^{\star}_{x_{i}}C$. 
\end {itemize}

\subsection {Plan of the paper.} The section following this introduction contains generalities about the spaces of stable maps. We collect there known results and we fix the notation. We also discuss the Hodge theory needed for our arguments. We will indicate the proof of Theorem $\ref {main}$ in the second section. The last section presents a few conjectures. The appendix contains a discussion of the higher genus tautological systems. 

\subsection {Acknowledgments.} We would like to thank Alina Marian for helpful conversations and enjoyable coffee breaks. We thank Professor Gang Tian for encouragement, support and guidance, Professor Aise Johan de Jong for several very useful discussions, and Professor Rahul Pandharipande for the interest shown in this work. 

\section {Generalities about the stable map spaces}

The purpose of this section is to collect various facts about the geometry of the stable map spaces which we will need later, and to fix the notation. We discuss the stratification with respect to the dual graphs and the associated Deligne spectral sequence. 

\subsection {The stratification by dual graphs.} \label{s1111}To begin, we let $X=G/P$ be a projective algebraic homogeneous space, where $G$ is a semisimple algebraic group and $P$ is a parabolic subgroup. We also fix $\beta\in H_{2}(X, \mathbb Z)$ a homology class. 

We let ${\overline M}_{0,n}(X, \beta)$ be the coarse moduli scheme of $n$ pointed stable maps to $X$ in the class $\beta$. A construction of the moduli scheme in algebraic geometry was given by Fulton and Pandharipande in \cite {FP}. It is shown there that the moduli scheme is a normal projective variety with finite quotient singularities - an orbifold if we work in the analytic category. In this subsection, we will discuss the coarse moduli schemes, but it should be clear how to extend our conclusions to the moduli stacks. 

To each stable map we associate the dual tree which carries degree labels and legs. We agree that a vertex $v$ has degree $\beta_v$ and $n(v)$ incident flags (i.e. edges and legs). The moduli space $\fscheme$ is a union of strata consisting of maps with fixed dual graph $\Gamma$. The closure of this stratum can be described as the image of a ramified covering of degree $Aut(\Gamma)$. The relevant morphism is: $$\zeta_{\Gamma}:\left(\prod_{v\in V(\Gamma)} \overline {M}_{0,n(v)}(X, \beta_v)\right)^{\text {Edge}(\Gamma)}\rightarrow \fscheme.$$ The left hand side is a fibered product along evaluation maps at the markings on the moduli spaces determined by the edges of $\Gamma$. 

There is a dense open stratum of maps with irreducible domains. Its complement is a union of divisors with normal crossings (up to a finite group action). However, this does not mean that the components of the boundary divisors do not self-intersect; in fact they almost always do. We emphasize this point for later reference when we write down the Deligne spectral sequence. 

The boundary strata are indexed by stable trees with one edge and two vertices. Each vertex has legs labeled by two sets $A, B$ with $A\cup B=\{1, \ldots, n\}$ and degrees $\beta_A, \beta_B$ adding up to $\beta$. Of course, stability means that if $\beta_A=0$ then $|A|\geq 2$, and similarly if $\beta_B=0$ then $|B|\geq 2$. The corresponding boundary stratum $$\iota:\overline {D}(\beta_A, \beta_B, A, B)\hookrightarrow \fscheme$$ is the image of a gluing map $$\zeta:{\overline M}(A,B, \beta_A, \beta_B)=\overline {M}_{0, {A \cup \{\bullet\}}}(X, \beta_A) {\times_{X}} \overline {M}_{0, \{\star\} \cup B}(X, \beta_B)\to \fscheme,$$ where $\bullet$ and $\star$ correspond to the double point of the domain curve. The following two cases may occur: \begin {itemize} \item If the tree has no automorphisms, then $\zeta$ is the normalization map of the boundary stratum $\overline {D}(A, B, \beta_A, \beta_B)$. If both $A, B\neq \emptyset$ then $\zeta$ is an isomorphism. Observe that if $A=\emptyset$ (or if $B=\emptyset$), then the corresponding divisor may self-intersect in a codimension two stratum of maps with three components glued at two nodes, the middle component containing all the marked points, the two external components having the same degree. 
\item If $A$ and $B$ are both empty and $\beta_A=\beta_B$, then the corresponding dual tree has a non-trivial automorphism. We need to factor out the $\mathbb Z/{2\mathbb Z}$ symmetry to get the normalization map of the corresponding boundary stratum.
\end {itemize}

\subsection {The Deligne spectral sequence} 

In this subsection, we review the main ingredients of Deligne's spectral sequence. Then, we apply the general theory to the case of the moduli space of stable maps. 

To get started, we let $Y$ be a smooth complex projective variety (or a projective orbifold later), $D$ be a divisor with normal crossings in $Y$; self-intersecting components are allowed. Let $U$ denote the complement of the divisor $D$ in $Y$, and let $j$ be the inclusion $U\hookrightarrow Y$. We denote by $D^p$ the subspace of $Y$ consisting of points of multiplicity at least $p$, and we let $\widetilde D_p$ be its normalization. Locally, $D$ is union of smooth divisors, and $D^p$ collects the points belonging to intersections of $p$ of them. We agree that $D^0=Y$. We will use cohomology with coefficients in the ``orientation" local system $\epsilon_p$. This local system is defined on $\widetilde D^{p}$ as follows. For points $y$ belonging to $p$ local components of $D$, we set $\epsilon_p$ to be the determinant of the space of the $p$ local components (see \cite {D} for a complete discussion). In the case when $D$ is union of smooth irreducible components without self-intersections, $\epsilon_p$ can be trivialized by choosing an ordering of the components. However, the case we will consider will involve self-intersecting components.

The cohomology of the open stratum $H^{\star}(U)$ carries a mixed Hodge structure which can be described explicitly in terms of the de Rham complex of logarithmic differentials. As a consequence of this general construction, $H^{n}(U)$ has a weight filtration $$0\subset W^n \subset W^{n+1} \subset \ldots \subset W^{2n}=H^n(U).$$ We will mainly be interested in the lowest piece of the filtration which can be computed as the restriction of the cohomology of $Y$: $$W^n=j^{\star} H^{n}(Y).$$ Additionally, there is another filtration whose role is to give the successive quotients $W^i/W^{i-1}$ a pure Hodge structure of weight $i$. 

There is a spectral sequence relating all the ingredients of the above discussion. Its $E_1$ term is $$ E_1^{-p,q}=H^{-2p+q}(\widetilde D^p, \epsilon_p),$$ and the first differentials of the spectral sequence are signed sums of Gysin inclusions. One of the main results of Hodge theory is that the higher differentials are all zero, and then the spectral sequence collapses to $E_{\infty}^{-p,q}=E_2^{-p,q}$ which is the piece of weight $q$ on $H^{q-p}(U)$ \cite {D}.

It was shown by Grothendieck \cite {Gr2} that the definition of the lowest Hodge piece is independent of the compactification $Y$ of $U$: for the purposes of defining $W^nH^n(U)$ we can pick any {\it smooth} compactification $Y$, maybe without normal crossings complement, and consider the restrictions of differential forms on $Y$. That is, for {\it all} smooth compactifications $j:U\hookrightarrow Y$, we have
$$W^{n}H^{n}(U)=j^{\star}H^{n}(Y).$$ It follows immediately that if $V\subset U$ have a common {\it smooth} compactification then the restriction map \begin{equation}\label{hsurj}W^n H^n(U)\to W^n H^n(V)\end{equation} is surjective. For example, if $V$ is Zariski open in a smooth variety $U$ such a compactification can always be found by Nagata's theorem and resolution of singularities. 

A similar remark can be proved in equivariant cohomology, if both $U$ and $V$ are equipped with compatible linearized actions of an algebraic group $G$. First, for any scheme $X$ with a $G$ action, the Hodge structure on the equivariant cohomology $H^{\star}_G(X)$ is constructed using the simplicial schemes $[X/G]_{\bullet}$ (see Chapter 6 in \cite {D} for notation). The equivariant version of Grothendieck's result is obtained as follows. It is proved in \cite {EG} that for each $n$, there exists an open subset $T$ of an affine $G$-space, equipped with a free $G$-action and whose complement has large codimension compared to $n$. The morphism of simplicial schemes $[X\times T/G]_{\bullet}\to [X/G]_{\bullet}$ induces a morphism of Hodge structures: $$H^{n}_{G}(X)=H^n([X/G]_{\bullet})\to H^n([X\times T/G]_{\bullet})= H^{n}_{G}(X\times T)=H^n(X\times_{G} T).$$ It is shown in \cite {EG} that the map above is an isomorphism. When $X$ is smooth, the morphism of {\it schemes} $X\times_{G} T\to T/G$ is also smooth \cite {EG}. For our applications, the base can always be chosen to be smooth, so $X\times_{G} T$ is also smooth. We combine the isomorphism of Hodge structures above with the surjectivity \eqref{hsurj} for the schemes $X\times_{G} T$ when $X$ is either $U$ and $V$. We conclude that: \begin {equation}\label {gro} W^{n}H^n_G(U)\to W^nH^n_G(V) \text { is surjective.} \end {equation} The similar statement about the Chow groups is evident. \\

We want to apply these general considerations to the space of stable maps. Although $\fscheme$ is not a smooth variety, its singularities are mild: they are all finite quotient singularities. There is an extension of Deligne's results to this setting which is worked out in \cite {St}. All the above results carry over without change. 

In the context of stable maps, we will mainly be concerned with the differential: $$d_1:E_1^{-1, k}=\bigoplus H^{k-2}\left(\overline {M}_{0, {A \cup \{\bullet\}}}(X, \beta_A) {\times_{X}} \overline {M}_{0, \{\star\} \cup B}(X, \beta_B)\right)^{-} \to E_1^{0,k}=H^{k}(\fscheme).$$ The cokernel of this map is the weight $k$ piece of $H^{k}(\fopenscheme)$ which we will proceed to identify in the next section when $X$ is a Grassmannian. 

The superscript ``-"on the cohomology groups of the boundary divisors comes from the orientation systems $\epsilon_p$. If the boundary graph has no automorphisms we consider the whole cohomology group. In the case of the $\mathbb Z/2\mathbb Z$ symmetry of the boundary graph we need to take fewer classes. In general, for any graph $\Gamma$, the boundary $\overline M(\Gamma)$ is dominated by a product of smaller moduli spaces of stable maps. The cohomology of the product carries a representation of $Aut(\Gamma)$. Each automorphism has a sign given by its action on the one dimensional space $det(\text {Edge}(\Gamma))$. The {\it minus} superscript indicates that we only look at classes which are anti-invariant under the sign representation of $Aut(\Gamma)$. For the boundary divisors of nodal maps with equal degrees on the branches, these are the $\mathbb Z/2\mathbb Z$-invariant classes.  

In particular, this discussion implies that the sequence: \begin {equation}\label {exactseq} \bigoplus H^{k-2}({\overline M}(A, B, \beta_A, \beta_B))\to H^{k}(\fscheme)\to W^k H^{k}(\fopenscheme)\to 0\end {equation} is exact (see also Corollary 8.2.8 in \cite {D}). The similar statement for the Chow groups is obvious. 

Once the exact sequence $\eqref{exactseq}$ is established, we can replace the coarse moduli schemes by the corresponding moduli stacks. To make sense out of the lowest piece of the Hodge structure on the smooth open stack $\mathcal {M}_{0,n}(X, \beta)$, we use the isomorphism with the cohomology of the coarse moduli scheme (Remark $1$ (iv)). Alternatively, the construction of a functorial mixed Hodge structure on the cohomology of algebraic stacks has been outlined in \cite {Dh}.

\section {Stable Maps to $SL$ flags}

In this section, we establish the proof of the main result. We begin with the case when the target space is a Grassmannian. We first identify the lowest piece of the Hodge structure on the cohomology of the open stratum of irreducible maps. We start with the case of three marked points, then move down to $0$, $1$ and $2$ markings. We conclude the argument by showing that the boundary classes are tautological. To this end, we prove a result about the cohomology of fibered products, essentially using techniques of \cite {D} and \cite {BF}. To complete the proof for general $SL$ flags, a discussion of the cohomology of the Hyper-Quot scheme is required.

\subsection {Stromme's description of the Quot scheme.} In the next few subsections we let $X$ be the Grassmannian of $r$ dimensional quotients of an $N$ dimensional vector space $V$. We will begin with a description of the smooth scheme of degree $d\geq 1$ morphisms to $X$, which we denote by $\text {Mor}_d(\p^1, X)$. When $X$ is the projective space, this discussion is trivial since the space in question can be described as an open subvariety in a projective space. However, for other flag varieties $X$ such a convenient description is not as easy to come by; it was obtained by Stromme for Grassmannians \cite {quot}. We also remark that when $X$ is a convex toric variety, the scheme of morphisms from $\p^1$ to $X$ admits a description similar to the one for $\pr$ \cite {Si}; it would be interesting to pursue our argument in this case. 

We will now explain Stromme's construction. To fix the notation, we let $\s$ and $\q$ denote the tautological subbundle and quotient bundle on the Grassmannian, sitting in the exact sequence: $$0\to \s\to V\otimes \mathcal O_{X}\to \q\to 0.$$  To give a degree $d$ morphism $f: \p^1 \to X$ is the same as giving a degree $d$, rank $r$ quotient vector bundle $F=f^{\star} \q$ as follows: $$V \otimes \mathcal O_{\p^1} \to F\to 0.$$ Allowing quotients which may not be locally free, we obtain the smooth compactification of $\mor$ which is known as Grothendieck's $Quot$ scheme. We will denote this by $Quot(N, r, d)$. It will be clear below that the cohomology of the $Quot$ scheme is easier to understand than that of the Kontsevich spaces. Incidentally, we note that this compactification was also used in \cite {BDW} to compute the ``Gromov invariants" of Grassmannians.

We consider two natural vector bundles $\mathcal A_{-1}$ and $\mathcal A_{0}$ on the scheme $\mor$. Their fibers over a morphism $f$ are $H^{0}(f^{\star} \q \otimes \o(-1))$ and $H^{0}(f^{\star}\q)$ respectively. These vector bundles extend over the $Quot$ scheme compactification. We let $\mathcal F$ be the universal quotient of the trivial bundle on $\p^1\times Quot$, and we let $\pi: \p^1\times Quot \to Quot$ be the projection on the second factor. The extensions are, for $m\geq -1$: $$\mathcal A_m=\pi_{\star}(\mathcal F\otimes \o(m)).$$ 

The relevance of these vector bundles for our discussion comes from the following consequence of Beilinson's spectral sequence. The bundle $F$ is the last term of the ``monad": 
\begin {center}
$\xymatrix{
0\ar[r] & H^{0}(F \otimes \o(-1))\otimes \o(-1)\ar[r] & H^{0}(F) \otimes  \o \ar[r] &F \ar[r] & 0.}
$
\end {center}

If we pick bases for $H^{0}(F(-1))$ and $H^{0}(F)$, we can identify these vector spaces with $W_{-1}=\c^d$ and $W_0=\c^{r+d}$ respectively. We rewrite the above diagram as follows: 
\begin {center}
$\xymatrix
{\;&\;& V\otimes \o\ar[d]\ar[dr]\\
0\ar[r] & W_{-1}\otimes \o (-1)\ar[r] & W_0 \otimes  \o \ar[r] &F \ar[r] & 0.}
$\end {center}

It is evident now that the datum of the quotient $V \otimes \o \to F$ can be encoded as an element of the affine space: $$H^{0}({\it Hom} (W_{-1}\otimes \o(-1), W_0\otimes \o))\oplus \text {Hom}(V, W_0)=P\otimes H^{0}(\o(1)) \oplus Q$$ where $P=\text {Hom} (W_{-1}, W_0)$ and $Q=\text {Hom} (V, W_0)$. Of course, not every element in this affine space is allowed; we need to impose the condition that the quotient we obtain be locally free and that the map from the trivial bundle $V\otimes \o\to F$ be surjective. In fact only an open subset of this affine space corresponds to morphisms $f\in\mor$. We also need to account for the $GL_d\times GL_{r+d}$ ambiguity coming from the action on the space of bases of $H^{0}(F(-1))\oplus H^{0}(F)$. \\

We carried out the above discussion on the level of closed points, but in fact Stromme's construction takes care of the scheme structure as well. We obtain:

\begin {fact}[Stromme, \cite {quot}] The total space $\mathcal T$ of the bundle of $GL_d\times GL_{r+d}$ frames of the vector bundle $\mathcal A_{-1} \oplus \mathcal A_{0}$ over $\mor$ sits as an open subscheme in the affine space $P\otimes H^{0}(\o(1)) \oplus Q$. In fact, this description extends to the $Quot$ scheme compactification of $\mor$. \\

As a consequence of equation $\eqref {gro}$ with trivial group, we conclude that the lowest piece of the Hodge filtration $W^{\star} H^{\star} (\mathcal T)=0.$
\end {fact}

To end the summary of Stromme's results, the following result of Grothendieck is needed. 

\begin {fact} [Grothendieck, \cite {Gr}] \label {gro2} Let $E$ be any rank $r$ vector bundle over a {\it smooth} compact base $X$. If $P$ denotes the bundle of $GL_r$ frames, then there is a surjective map: $$ H^{\star}(X) \to W^{\star}H^{\star}(P)$$ whose kernel is the ideal generated by the Chern classes $c_i(E)$. The same statement holds equivarianly for an algebraic group action, and for the Chow rings.
\end {fact}

The statement in \cite {Gr} is written in Chow, but the same proof also works in cohomology. The strategy is to reduce to the case of line bundles. The equivariant case can be deduced via a similar reasoning to that used to prove equation \eqref{gro}. 

Applying these two facts to our setting we obtain that $H^{\star} (Quot)$ is generated by the Chern classes $c_i(\mathcal A_{0})$ and $c_i(\mathcal A_{-1})$. Therefore, using $\eqref {hsurj}$, we obtain:

\begin {fact} [Stromme, \cite {quot}] The lowest weight Hodge piece $W^{\star}H^{\star}(\mor)$ is generated by the Chern classes $c_i(\mathcal A_0)$ and $c_i(\mathcal A_{-1})$. The same result holds for the Chow rings.
\end {fact}
\subsection {At least three marked points.} Our goal is to prove the following statement about the open stratum $\mathcal M_{0,n}(X, d)$: 

\begin {lemma}\label {nomarks}
The lowest piece of the Hodge structure $W^{\star}H^{\star}({\mathcal M}_{0,n}(X,d))$ is spanned by restrictions of the tautological classes on $\overline {\mathcal M}_{0,n}(X, d)$ of Definition $1$. The same result holds for the Chow groups of ${\mathcal M}_{0,n}(X, d)$. 
\end {lemma}

In this subsection, we will consider the case of the above lemma when $n\geq 3$. We use the results proved above about the $Quot$ scheme. It will be useful to consider the sheaves on the stack $\fstack$ defined, for $m\geq 0$, as: $$\mathcal G_m=\pi_{\star}ev^{\star} \left(\q\otimes (\det \q)^{\otimes m}\right).$$ Here $\pi$ and $ev$ are the projection from the universal curve and the evaluation map respectively. We will be interested in the Chern classes of the restrictions $j^{\star} \mathcal G_m$ to the open stratum $j:\mathcal M_{0,n}(X, d)\to \overline {\mathcal M}_{0,n}(X, d)$. We show that these generate the lowest Hodge piece of cohomology (or the Chow groups). 

We consider the fiber diagram, where $\mathcal C_{0,n}(X, d)$ is the universal curve over $\mathcal M_{0,n}(X, d)$:
\begin {center}$
\xymatrix{
\mathcal C_{0,n}(X, d) \ar[r]^{\hat q}\ar[d]^{\pi}\ar@/^1pc/[rr]|{ev}& \p^1\times \mor \ar[r]^{\hspace{.4in} ev}\ar[d]^{\hat \pi}& X\\
\mathcal M_{0,n}(X, d) \ar[r]^{q}& \mor}$
\end {center}
In the above diagram, $q:{\mathcal M}_{0,n}(X, d)\to \mor$ is the forgetful map defined by forgetting all but the first three points, and $\hat q$ is defined similarly. 
We let $$p:\p^1\times \mor \to \p^1$$ be the projection. Since $ev^{\star} \det \q$ and $p^{\star}\o(d)$ agree on the fibers of $\hat \pi$, there exists a line bundle $\mathcal L$ on $\mor$ for which the following equation is satisfied: $$ev^{\star} \det \q = p^{\star}\o(d)\otimes \hat \pi^{\star} \mathcal L.$$ We compute: \begin {equation}\label {jstar}
j^{\star} \mathcal G_m=\pi_{\star}ev^{\star}\left(\q\otimes (\det \q)^{\otimes m}\right)=\pi_{\star} \hat q^{\star} \left(ev^{\star} \q \otimes p^{\star} \o(dm) \otimes \hat \pi^{\star} \mathcal L^{\otimes m}\right)=\end {equation}
$$=q^{\star} \hat  \pi_{\star}\left( ev^{\star} \q \otimes p^{\star} \o(dm)\otimes \hat \pi^{\star} \mathcal L^{\otimes m}\right)=q^{\star} \mathcal A_{dm} \otimes q^{\star}\mathcal L^{\otimes m}.$$

For $n\geq 3$, there is a bijective morphism $\fopenscheme\to M_{0,n}\times \mor$ induced by the forgetful map to $M_{0, n}$ and the morphism $q$. Since the domain and target are smooth, we obtain an isomorphism. Kunneth decomposition can be used to understand the Hodge structure on the open stratum: $$H^{k}({\mathcal M}_{0,n}(X, \beta))=\bigoplus_{i+j=k} H^{i}(M_{0,n}) \otimes H^{j}(\mor).$$ In \cite {zero}, it is proved that the $i^{\text {th}}$ cohomology $H^{i}(M_{0,n})$ carries a pure Hodge structure of weight $2i$. Since $\mor$ is smooth, its $j^{\text {th}}$ cohomology carries weights between $j$ and $2j$. Hence to get the weight $k$ piece $W^kH^{k}(\mathcal M_{0,n}(X, \beta))$ we need $i=0$, $j=k$. Thus, this weight $k$ piece is isomorphic to the the weight $k$ piece $W^kH^{k}(\mor)$ via the pullback: $$q^{\star}:W^{\star} H^{\star} (\mor)\to W^{\star} H^{\star}(\mathcal M_{0,n}(X, d)).$$ 

It follows from {\it Fact $3$} and the exact sequence \cite {quot}: $$0\to\mathcal A_k \to \mathcal A_{k+1}^{\oplus 2}\to \mathcal A_{k+2} \to 0, \;\; \text { for } k\geq 0,$$ that $W^{\star}H^{\star}(\mor)$ is generated by the Chern classes of $\mathcal A_{0}$ and $\mathcal A_{dm}$, for any $m\geq 1$. Better, we can pick as generators the Chern classes of $\mathcal A_{dm}\otimes \mathcal L^{\otimes m}$ for all values $m\geq 0$. To see this, observe that we can express the Chern classes of the bundles $\mathcal A_{dm}$ in terms of those of $\mathcal A_{dm}\otimes \mathcal L^{\otimes m}$. Indeed, it suffices to prove that $c_1(\mathcal L)$ can be expressed in terms of the collection of Chern classes $c_1(\mathcal A_{dm}\otimes \mathcal L^{\otimes m})$ for all $m\geq 0$. This follows from the Chern class computation: $$c_1(\mathcal A_{dm}\otimes \mathcal L^{\otimes m})=mc_1(\mathcal L) \text { rank }{\mathcal A}_{dm}+c_1(\mathcal A_{dm})=m^2 d r c_1(\mathcal L)+\text {linear terms in } m.$$ The isomorphism $q^{\star}$ above and equation $\eqref {jstar}$ show that the Chern classes of $j^{\star} \mathcal G_m$ generate $W^{\star} H^{\star}(\mathcal M_{0,n}(X, d))$ for $n\geq 3$. 

Finally, to completely prove Lemma $\ref {nomarks}$ when $n\geq 3$, we need to show that the Chern classes of $\mathcal G_m$ on $\mathcal M_{0,n}(X, d)$ are restrictions of the tautological classes of Definition $1$. We repeat the argument of Lemma $11$ in \cite {O1}. There, we explained the required Mumford-Grothendieck-Riemann-Roch computation. To conclude, we need to observe that the Chern class of the relative dualizing sheaf $c_1(\omega_{\pi})$ is tautological. Here $\pi$ is the forgetful morphism. This is a consequence of Proposition $1$ in \cite {O2} which shows all codimension $1$ classes on $\overline {\mathcal M}_{0, n}(X, d)$ are tautological. One can argue differently by using the Pl\"ucker embedding to reduce the statement to the case of a projective space $\p^N$ in which the Grassmannian embeds. Then, we invoke \cite {divisors} which gives expressions for $c_1(\omega_{\pi})$ in terms of tautological classes as in Definition 1. Restricting to $X$ yields the claim. 

The statement for the Chow groups follows in the same fashion. We need the observation that for any scheme $T$, the map: $$q^{\star}:A_{\star}(T)\to A_{\star}(\mon\times T)$$ is an isomorphism. This is well known. For example, it follows from comparing the exact sequences:

{\xymatrix{
\bigoplus_i A_{\star}(\overline D_i) \otimes A_{\star}(T)\ar[r]\ar[d] &A_{\star}(\monbar)\otimes A_{\star} (T)\ar[r]\ar[d] &A_{\star}(\mon)\otimes A_{\star}(T) \cong A_{\star} (T)\ar[r]\ar[d]^{q^{\star}}& 0\\
\bigoplus_i A_{\star}(\overline D_i \times T)\ar[r] & A_{\star}(\monbar \times T)\ar[r]& A_{\star}(\mon\times T)\ar[r]& 0}}
\noindent Here $\overline D_i$ are the boundary divisors of $\monbar$, which are products of lower dimensional moduli spaces of stable marked rational curves. Since the first two vertical arrows are isomorphisms as it is shown in Section $2$ of \cite {K}, the third vertical arrow is also an isomorphism.

\subsection {Fewer marked points.} We will now prove Lemma $\ref{nomarks}$ when the domain has fewer marked points. A case by case analysis depending on the number of markings is required.   

\subsubsection {No marked points.} To begin, we consider the case of no marked points. We let $SL_2$ act on $\c^2$ in the usual way. In turn, we obtain an action on the scheme of morphisms $\mor$: $$SL_2\times \mor \ni (g, f)\to f\circ g^{-1}\in \mor.$$ Since each morphism is finite onto its image, it is easy to derive that the action has finite stabilizers. The $PSL_2$ quotient of $\mor$ equals the coarse space $M_{0,0}(X, d)$ and since the center in $SL_2$ acts trivially, we see that the space of $SL_2$ orbits of $\mor$ is the topological space underlying $M_{0,0}(X,d)$. A well known result, which is proved for example in \cite {Br}, gives an isomorphism between the cohomology of the orbit space and the equivariant cohomology. In our case, this translates into an isomorphism:\begin {equation}\label {isom} H^{\star}(M_{0,0}(X, d))=H^{\star}_{SL_2}(\mor).\end {equation} Even more, the right hand side can be given a Hodge structure using simplicial schemes \cite {D}, Chapter $6$, which, by functoriality is compatible with the structure on the left hand side. 

We now move the discussion to the algebraic category. It is easy to see using the numerical criterion of stability that the action of $SL_2$ on $$\text {Mor}_d(\p^1, \p^N)\hookrightarrow \p^{Nd+N+d}$$ has only stable points. Indeed, fixing a $1$-parameter subgroup $\lambda(t)=\text{diag}(t^{a}, t^{-a})$, and letting $$f=[f_0: \ldots: f_N], \quad f_k=\sum_{i=0}^{d} c_{i}^{(k)}z^i w^{d-i}$$ be the components of a degree $d$ morphism, we need to show $$\min\{2i-d: c_{i}^{(k)}\neq 0 \text{ for some }k\}<0<\max\{2i-d: c_{i}^{(k)}\neq 0\text{ for some }k\}.$$ This is clear since $f_0, \ldots, f_N$ have no common roots. The same statement then holds for $\mor$ using the Pl\"ucker embedding. In turn, this implies the existence of a geometric quotient $\mor/SL_2$. This quotient can be identified with $M_{0,0}(X, d)$. Indeed, the natural map $$\text{Mor}(\mathbb P^1, X)\to M_{0,0}(X, d)$$ is $SL_2$-invariant, so we obtain a morphism $\text{Mor}(\mathbb P^1, X)/SL_2\to M_{0,0}(X, d)$ which is bijective on points. Since $M_{0,0}(X, d)$ is normal \cite{FP}, the above is an isomorphism. 

We will make use of the isomorphism \cite {EG}: \begin {equation} \label {isom1} A^{k}(M_{0,0}(X, d))=A^{k}_{SL_2}(\mor)=A^{k}(\mor\times_{SL_2} W).\end {equation} In the topological category, we compute equivariant cohomology taking for $W$ the contractible space $ESL_2$. In the algebraic case, we let $W=W_k$ be any a smooth open subvariety of an affine $SL_2$-space, which has large codimension compared to $k$ and which has a free $SL_2$ action (cf. \cite {EG}, see also the previous section). 

Let us write, for now, $G=SL_2$; later, we will make use of other groups as well. We will use the following standard notation. For any $G$ scheme $X$, $X_G$ will denote the equivariant Borel construction. In the algebraic setting, $X_G$ will stand for any of the mixed {\it schemes} $X\times_{G} W$ (cf. \cite {EG}) where $W$ is described in the previous paragraph. We observed in the first section that all equivariant models $X_G$ can be chosen to be smooth.  Generally, properties of equivariant morphisms such as smoothness, or flatness and properness, still hold for the induced morphisms between the mixed spaces \cite {EG}; we will make use of this fact below. Moreover, any $G$-linearized bundle $E\to X$ lifts to a bundle $E_G\to X_G$ whose total space is $E\times_{G} W$. 

We will compute the right hand side of $\eqref {isom}$ and $\eqref {isom1}$ using arguments similar to Stromme's. We will need to extend the $SL_2$ action to the $Quot$ scheme and to lift it to the bundles $\mathcal A_{-1}, \mathcal A_0$. For example, the $SL_2$ linearization of $\mathcal A_{-1}$ is determined by the usual $SL_2$ linearization of $L=\o(-1)$. To be precise, we consider the following diagram:
\begin{center}$\xymatrix {
SL_2 \times \p^1\ar@<.5ex>[r]^{\sigma}\ar@<-.5ex>[r]_{\alpha}& \p^1\\
SL_2 \times \p^1 \times Quot \ar[d]^{\pi} \ar[u]^{p} \ar@<.5ex>[r]^{\;\;\sigma}\ar@<-.5ex>[r]_{\;\;\alpha} & \p^1 \times Quot \ar[u]_{p}\ar[d]^{\pi}\\
 SL_2 \times Quot \ar@<.5ex>[r]^{\sigma}\ar@<-.5ex>[r]_{\alpha}& Quot.
}$\end {center}
Here, $\sigma$ is the action, while the morphisms $\alpha$, $\pi$ and $p$ are projections. We have an isomorphism $\eta:\sigma^{\star} \o(-1)\to \alpha^{\star} \o(-1)$. Writing as before $\mathcal F$ for the universal quotient sheaf on $\p^1 \times Quot$, the linearization of $\mathcal A_{-1}$ is obtained as follows: 
\begin {eqnarray*} \sigma^{\star} {\mathcal A}_{-1}&=&\sigma^{\star} \pi_{\star} (\f \otimes p^{\star} \o(-1))=\pi_{\star} \sigma^{\star} (\f\otimes p^{\star} \o(-1))\\
&\cong& \pi_{\star} (\alpha^{\star} \f\otimes p^{\star} \sigma^{\star} \o(-1)) \cong \pi_{\star} (\alpha^{\star} \f \otimes p^{\star} \alpha^{\star} \o(-1))\\ &=& \pi_{\star} \alpha^{\star}(\f \otimes p^{\star} \o(-1))=\alpha^{\star} \pi_{\star}(\mathcal F \otimes p^{\star} \o(-1))=\alpha^{\star} \mathcal A_{-1}.
\end {eqnarray*}

We will be concerned with the mixed space $Quot_{SL_2}$ and the bundles $\mathcal A_m^{SL_2}=\mathcal A_m\times_{SL_2}W$ which are the lifts of the equivariant bundles $\mathcal A_m$, for $m\in\{-1, 0\}$. A moment's thought shows that the Stromme embedding described in {\it {Fact $1$}} is $SL_2$ equivariant. It is immediate that taking frames of a vector bundle commutes with the construction of the mixed spaces and of mixed bundles over them. Then, the bundle of (split) frames of $\mathcal A_{-1}^{SL_2}\oplus \mathcal A_0^{SL_2}$ is the mixed space $\mathcal T_{SL_2}$ and moreover, it can be realized as a subscheme \begin{equation}\label{subsch}\mathcal T_{SL_2} \hookrightarrow (P\otimes H^{0}(\p^1, \o(1)))\oplus Q)_{SL_2}.\end{equation} 

The latter space can be described explicitly. In the topological category, $BGL_2$ can be realized as the infinite Grassmannian $\mathbb G$ of $2$ dimensional planes endowed with a tautological rank $2$ bundle $\mathbb S$. Its frame bundle serves as a model for $EGL_2=ESL_2$. In the algebraic category, we consider the truncated models of the infinite dimensional constructions. For instance, $W/GL_2$ will be a finite dimensional Grassmannian of $2$ dimensional planes. Then, $W$ will be the bundle of $GL_2$ frames of the tautological rank $2$ bundle $\ms$ over $W/GL_2$. $W/SL_2$ is the bundle of $\cs$ frames of the determinant $\det \;\ms\to W/GL_2$. This is a consequence of the fact that any frame of $\ms$, say $w\in W$, gives rise to a frame $\det w$ of $\det \ms$. We will denote by $\tilde \ms$ the pullback to $W/SL_2$ of the tautological bundle $\ms$. 

It is then clear that the second mixed space in \eqref{subsch}, which we are trying to describe, can be realized as the total space of the bundle: $$P\otimes \tilde {\mathbb S}^{\star}\oplus \mathcal O\otimes Q\to W/SL_2.$$ Therefore, its cohomology (or Chow rings) can be computed from that of the base. In turn, this is easily seen to be (via {\it Fact $2$}): $$W^{\star}H^{\star}_{SL_2}=H^{\star}_{GL_2}/c_1(\Lambda^{2}\tilde{\mathbb S})=\mathbb C\left[c_2(\tilde{\mathbb S})\right]$$ 

Repeating the arguments of the previous subsection in the equivariant setting (making use of {\it Fact} $2$ and equation $\eqref {gro}$) we obtain two surjections: $$W^{\star}H^{\star}_{SL_2}\to W^{\star}H^{\star}_{SL_2}(\mathcal T),\;\; H^{\star}_{SL_2}(Quot)\to W^{\star}H^{\star}_{SL_2}(\mathcal T).$$ We have explicit generators for the kernel of the second map, namely the equivariant Chern classes of the bundles $\mathcal A_{-1}, \mathcal A_{0}$. We conclude the following:\vskip.1in
\begin{claim} The cohomology $H^{\star}_{SL_2}(Quot)$ is generated by the equivariant classes of $\mathcal A_0$ and $\mathcal A_{-1}$ together with the pullback of the Chern class $4 c_2(\tilde{\mathbb S})=c_2(Sym^{2} \tilde{\mathbb S})$ from $W/SL_2$. Therefore, the same is true about the lowest piece of the Hodge structure of $H^{\star}_{SL_2}(\mor)$. \end{claim}
We are now ready for the proof of Lemma $\ref {nomarks}$. We consider the following diagram of fiber squares:
\begin {center}
$\xymatrix{
\p(\tilde {\ms}) = \p^1\times_{SL_2} W\ar[d]^{\hat \pi} &({\p^1\times \mor})\times_{SL_2} W \ar[l]_{\hspace{-.3in}\hat q}\ar[r]^{\hspace{.6in} \tilde \epsilon}\ar[d]^{\tilde\pi} & \mathcal C \ar[r]^{ev}\ar[d]^{\pi}& X\\
W/SL_2&{\mor\times_{SL_2} W}\ar[l]_{\hspace {-.3in}q}\ar[r]^{\hspace{.6in}\epsilon} & \mathcal {M}_{0,0}(X,d)\ar[r]^{p} & M_{0,0}(X, d)
}$
\end {center}

In the above diagram, we start with the flat algebraic family $\tilde \pi$ whose fibers are irreducible genus $0$ curves (which are not canonically identified to $\p^1$ because we factored out the $SL_2$-action). The classifying map to the moduli stack is denoted by $\epsilon$. In the same diagram, $\mathcal C$ is the universal curve, and $p:\mathcal M_{0,0}(X, d)\to M_{0,0}(X, d)$ is the natural morphism to the coarse moduli scheme.  

We know that $p$ induces isomorphisms in rational cohomology \cite {Be}. As noted in Section \ref{s1111}, the same is true about $p\epsilon$ in each fixed degree, provided $W$ is the appropriate truncation of the infinite dimensional model. Thus the conclusion also holds for $\epsilon$. It suffices to show that any class of lowest Hodge weight on the smooth space $\mor\times_{SL_2} W$ is the pullback from $\mathcal M_{0,0} (X, d)$ of restrictions of the tautological classes of Definition $1$.  

We will use the conclusions summarized in the {\it Claim} above. It is clear that $\mathcal A_0^{SL_2}$ corresponds to the bundle $\mathcal G_0=\pi_{\star}(ev^{\star}\q)$ under the isomorphism $\eqref {isom}$ induced by $\epsilon$. 

We argue that the Chern classes of the bundle $\mathcal A_{-1}^{SL_2}$ also come as pullbacks under $\epsilon$ of tautological classes on $\mathcal M_{0,0}(X,d)$. The argument for $\mathcal A_0^{SL_2}$ is similar. Let $L$ be the lift of the linearized bundle $p^{\star} \o(-1)$ on $\p^1\times \mor$ to the equivariant mixed space $(\p^1\times \mor)\times_{SL_2}W$. 
First we note that $\mathcal A_{-1}^{SL_2}=\tilde \pi_{\star} (\tilde \epsilon^{\star}ev^{\star}\q\otimes L).$ 

To compute the Chern classes of $\mathcal A_{-1}^{SL_2}$ we use Grothendieck-Riemann-Roch. The Chern classes are determined by the Chern character, which in turn is given by $$\tilde \pi_{\star}\left(\tilde \epsilon^{\star} ev^{\star} ch(\mathcal Q)\cdot e^{c_1(L)}\cdot \text{Todd}({\tilde \pi})\right).$$ Note that $L=\widehat q^{\star} \widehat L$ where $\widehat L$ is the analogous bundle over $\mathbb P^1\times_{SL_2} W.$ We note below that $\omega_{\widehat \pi}=\widehat L^{\otimes 2}$. Thus $$c_1(L)=\frac{1}{2} \widehat q^{\star} c_1(\omega_{\widehat \pi})=\frac{1}{2} \tilde \epsilon^{\star} c_1(\omega_{\pi}).$$ Writing $$\tau=ev^{\star}ch(\mathcal Q)\cdot e^{\frac{1}{2}c_1(\omega_{\pi})}\cdot \text{Todd} ({\pi}),$$ the above Chern character becomes $$\tilde \pi_{\star}\tilde \epsilon^{\star} \tau =\epsilon^{\star} \pi_{\star} \tau.$$ The class $\tau$ and the pushforward $\pi_{\star}\tau$ are tautological since the class $c_1(\omega_{\pi})$ is so, as we noted in Section 2.2.

Finally, we discuss $c_2(Sym^2 \tilde{\mathbb S})$. We claim that on $\mor\times_{SL_2} W$ the following equation holds true: \begin {equation}\label {sym}\epsilon^{\star} \pi_{\star} \omega^{\star}_{\pi}=q^{\star}Sym^{2}\tilde{\mathbb S}^{\star}.\end {equation} 
The Chern classes of $\pi_{\star} \omega_{\pi}^{\star}$ are tautological by the usual argument involving Grothendieck-Riemann-Roch. 

To establish the equation above, observe that we are in the particularly favorable situation when all relative dualizing sheaves involved are line bundles over smooth bases, hence taking duals causes no problems. We have $\epsilon^{\star} \pi_{\star} \omega^{\star}_{\pi}=\tilde \pi_{\star}\omega^{\star}_{\tilde \pi}=q^{\star} \hat \pi_{\star}\omega^{\star}_{\hat \pi}$. The following Euler sequence on $\p (\tilde{\mathbb S})$: $$0\to \mathcal O \to \hat \pi^{\star} \tilde{\ms} \otimes \mathcal O_{\p(\tilde{\ms})}(1)\to \omega_{\hat \pi}^{\star}\to 0$$ shows that  $$\omega^{\star}_{\hat \pi}= \det (\hat \pi ^{\star}{\tilde\ms}\otimes \mathcal O_{\p(\tilde{\ms})}(1))=\mathcal O_{\p(\tilde{\ms})}(2).$$ Here we used that $W/SL_2$ is the space of frames for $\Lambda^2 \ms$, so the pullback of $\Lambda^2 \ms\to W/GL_2$ to $W/SL_2$ is trivial. We immediately obtain $\hat\pi_{\star}\omega^{\star}_{\hat\pi}=Sym^{2} \tilde{\ms}^{\star}$, establishing $\eqref {sym}$. 

\subsubsection {One marked point.} The remaining two cases are entirely similar. For the case of one marking, we will use the action of the subgroup $N$ of $SL_2$ of matrices: $$N=\left\{\left(\begin{array}{cc} a & b \\ 0 & a^{-1} \end {array}\right)\big{|}\;\; a \in \cs, b\in \c\right\}.$$ We carry out the equivariant arguments of the previous section replacing $SL_2$ with the group $N$. 

We will identify $EN$ and $BN$. To construct algebraic families, we will work with the finite dimensional approximations of the topological models. As before, $W/GL_2$ will be a Grassmannian of $2$ dimensional planes of large, but finite, dimension. We identify $W$ with the bundle of $GL_2$ frames of the tautological bundle $\ms$. We have seen that $W/SL_2$ is the bundle of $\cs$ frames in $\det \; \ms.$ Moreover, the space $W/SL_2$ comes equipped with the tautological pullback bundle, which we denote by $\tilde \ms$. The space $W$ can be used to compute the $N$-equivariant cohomology. Then, the quotient $W/N=W\times_{SL_2} SL_2/N$ exists, and it will be the projective bundle $\mathbb P(\tilde \ms)$ over $W/SL_2$. In addition, letting $\tau: W/N\to W/GL_2$ be the projection, we have $\p^1\times_{N} W = \p(\tau^{\star} \ms)$. The last two statements follow by noting the map: $$\{\text {frames in a $2$ dimensional vector space $S$}\}/N \to \mathbb P^{1}(S),\; \text {frame \{$e,f$}\}\to \text {line spanned by $e$}.$$
The relevant diagram of fiber squares is:
\begin {center}
$\xymatrix{
\p(\tau^{\star}\ms)=\p^1\times_{N} W\ar[d]^{\hat\pi} &({\p^1\times \mor})\times_{N} W \ar[l]_{\hspace{-.2in}\hat q}\ar[r]^{\hspace{.5in}\tilde \epsilon}\ar[d]^{\tilde \pi} & \mathcal C \ar[r]^{ev}\ar[d]^{\pi}& X\\
\p(\tilde \ms)=W/N\ar@<1ex>[u]^{\hat z}&{\mor}\times_{N} W\ar[l]_{\hspace{-.2in}q}\ar[r]^{\hspace{.5in}\epsilon}\ar@<1ex>[u]^{\tilde z} & \mathcal M_{0,1}(X,d)\ar@<1ex>[u]^{z}\ar[r]^{p}& M_{0,1}(X, d)
}$
\end {center}

We will denote by $\xi$ the dual hyperplane bundle of $\p(\tilde \ms)\to W/SL_2$ and similarly $\Xi$ will be the dual hyperplane bundle of $\p(\tau^{\star}\ms)\to W/N$. The morphism $\hat \pi: \p(\tau^{\star}\ms)\to \p(\tilde \ms)$ has a canonical section $\hat z$ such that $\hat z^{\star}\Xi=\xi$. The flat algebraic family $\tilde \pi$ also has a section $\tilde z$ whose image in each fiber is the $N$-invariant basepoint $[1:0]\in \p^1$. We let $\epsilon$ be the classifying map to the moduli stack $\mathcal M_{0,1}(X, d)$.

The above description of $W/N$ as a projective bundle over $W/SL_2$ shows that $$W^{\star}H^{\star}_N=W^{\star}H^{\star}_{SL_2}[c_1(\xi)]/(c_1(\xi)^2+c_2(\tilde \ms))=\c[c_1(\xi), c_2(\tilde \ms)]/(c_1(\xi)^2+c_2(\tilde \ms))=\c[c_1(\xi)].$$  To complete the proof of Lemma $\ref {nomarks}$, we use the arguments in the previous subsection. We only have to write the class $q^{\star}c_1(\xi)$ on $\mor\times_N W$ as the pullback of a tautological class on $\mathcal M_{0,1}(X, d)$. 

The Euler sequence on $\p(\tau^{\star}\ms)$ shows that $$0\to \mathcal O\to \hat \pi^{\star}\tau^{\star} \ms \otimes \Xi \to \omega_{\hat \pi}^{\star}\to 0.$$ Taking determinants and recalling that $\tau^{\star}\det \ms$ is trivial, we obtain $\omega^{\star}_{\hat \pi}=\Xi^{\otimes 2}$. Therefore, $$\hat z^{\star} c_1(\omega_{\hat \pi})=-2c_1(\hat z^{\star} \Xi)=-2c_1(\xi).$$ It is then clear that: $$-2q^{\star} c_1(\xi)=q^{\star} \hat z^{\star} c_1(\omega_{\hat \pi})=\tilde z^{\star} c_1(\omega_{\tilde \pi})=\epsilon^{\star} (z^{\star} c_1(\omega_{\pi}))=\epsilon^{\star} \psi_1.$$ The proof is now complete, since the $\psi$ classes are tautological in genus $0$ by Lemma $2.2.2$ in \cite {divisors}. The latter expresses the $\psi$'s in terms of evaluation classes, boundaries and the $\kappa$ classes, in the case of projective spaces. From here the statement for Grassmannians $X$ follows as well by pullback under the Pl\"ucker embedding. 

\subsubsection {Two marked points.} The argument is again similar to the case of one marked point. We let $\cs$ act on $\p^1$ as follows $t\cdot [z:w]=[t^{-1}z:t w]$. We obtain the isomorphism: $$H^{\star}(\mathcal M_{0,2}(X, d))=H^{\star}_{\cs} (\mor).$$ We compute the right hand side. Here $W/\cs$ will be a large projective space and $W$ the bundle of $\cs$ frames of the tautological line bundle $\ms\to W/\cs$. 

The following diagram of fiber squares will be useful in our computation: 

\begin {center}
$\xymatrix{
\p(\ms\oplus \ms^{\star}) = \p^1\times_{\cs} W\ar[d]^{\hat\pi} &({\p^1\times \mor})\times_{\cs} W \ar[l]_{\hat q}\ar[r]^{\hspace{.5in}\tilde \epsilon}\ar[d]^{\tilde \pi} & \mathcal C \ar[r]^{ev}\ar[d]^{\pi}& X\\
W/\cs \ar@<.8ex>[u]\ar@<1.6ex>[u]&{\mor}\times_{\cs} W\ar[l]_{q}\ar@<.8ex>[u]\ar@<1.6ex>[u]\ar[r]^{\hspace{.5in}\epsilon} & \mathcal M_{0,2}(X,d)\ar@<.8ex>[u]\ar@<1.6ex>[u]\ar[r]^{p}& M_{0,2}(X, d)
}$
\end {center}
The family $\hat \pi$ has two tautological sections $\hat z, \hat w$ such that $\hat z^{\star} \xi=\ms$ and $\hat w^{\star} \xi=\ms^{\star}$; here, $\xi$ is the dual hyperplane bundle of $\p(\ms \oplus \ms^{\star})\to W/\cs$. The family $\tilde \pi$ also has two sections, their images are the two $\cs$ invariant points of $\p^1$ in each fiber.

We obtain generators for the lowest piece of the Hodge structure on $H^{\star}_{\cs}(\mor)$. In the light of the previous discussion, we will only need to explain that the class $q^{\star} c_1(\ms)$ coming from $W/\cs$ is a pullback of a tautological class on $\mathcal M_{0,2}(X, d)$. The argument is identical to the one in the previous section. From the Euler sequence along the fibers of $\hat \pi$, it follows that $c_1(\omega_{\hat \pi})=-2c_1(\xi)$. The computation below finishes the proof: $$-2 q^{\star} c_1(\ms)=-2q^{\star} \hat z^{\star}c_1(\xi)= \tilde z^{\star} \hat q^{\star} c_1(\omega_{\hat \pi})=\tilde z^{\star}c_1(\omega_{\tilde \pi})=\epsilon^{\star} z^{\star} c_1(\omega_{\pi})=\epsilon^{\star} \psi_1.$$

\subsection {General SL flag varieties.} Let us now consider the case of a general $SL$ flag variety $X$  parameterizing $l$ successive quotients of dimensions $r_1, \ldots, r_l$ of some $N$ dimensional vector space $V$. 

Pulling back the tautological sequence on $X$
\begin {equation}\label {taut}
0\to \s_1\to \ldots \to \s_l\to V\otimes \mathcal O_{X}\to \q_1\to \ldots \to \q_l\to 0.
\end {equation}
under a morphism $f:\p^1\to X$, we obtain a sequence of locally free quotients $$V\otimes \o \to F_1\to \ldots \to F_l\to 0$$ of ranks $r_i$ and degrees $d_i$. Allowing for arbitrary (not necessarily locally free) quotients, we obtain the Hyper-Quot scheme compactification $HQuot$ of the space of morphisms $\text {Mor}_{\beta}(\p^1, X)$. We will use the more explicit notation $HQuot (N, {\bf r}, {\bf d})$ when necessary. Here ${\bf r}=(r_1, \ldots, r_l)$ and ${\bf d}=(d_1, \ldots, d_l)$.

Since $HQuot$ is a fine moduli scheme, there is a universal sequence on $\p^1\times HQuot$: $$0\to \mathcal E_1\to\ldots \to \mathcal E_l \to V\otimes \mathcal  O \to \mathcal F_1\to \ldots \to \mathcal F_l\to 0.$$ We seek to show that:
\begin {lemma} \label {hq}The cohomology (and the Chow rings) of $HQuot (N, {\bf r}, {\bf d})$ is generated by the Kunneth components of $c_j(\mathcal F_i)$. 
\end {lemma}

Unfortunately, the arguments of the previous subsections do not extend to the present case. Even though a description of the Hyper-Quot scheme similar to the one in Section $2.1$ does exist \cite {Kim}, nonetheless we obtain an embedding of a principal bundle over $HQuot$ into a $\it {singular}$ affine variety. The existence of singularities is a serious (and not the only) obstacle in extending the proofs in Section $2.2$ to our new setting, as it is harder to keep track of the Hodge structures.

Nevertheless, the proof of the lemma stated above should be well known, but we could not find a suitable reference. To prove it, we will use a well-known trick of Beauville-Ellingsrud-Stromme \cite {ES}. Some of the details appear in the next subsection. There are difficulties in applying the same argument equivariantly. Instead, we will use a combination of the Leray spectral sequence and the already established equivariant results for the $Quot$ scheme.  We prove the following:

\addtocounter {lemma}{-2}
\begin{lemma}{\bf (bis)} \label {nomarksflags}
For any flag variety $X$ and any degree $\beta$, the lowest piece of the Hodge structure $W^{\star}H^{\star}(\mathcal M_{0,n}(X,\beta))$ is spanned by the restrictions of the tautological classes on $\overline {\mathcal M}_{0,n}(X, \beta)$ of Definition $1$. The same results hold for the Chow groups. 
\end {lemma}

To see this, we will assume Lemma $\ref {hq}$. We consider the following product of forgetful morphisms: $$i=\prod_{j}i_j: HQuot_{\p^1}(N, {\bf r}, {\bf d})\to \prod_{j} Quot(N, r_j, d_j).$$ We put together Lemma $\ref {hq}$, and the observation that the universal quotients on $\p^1\times Quot(N, r_j, d_j)$ pull back to the universal sheaves $\mathcal F_j$ on $\p^1\times HQuot$. We conclude that the pullback map $$i^{\star}:H^{\star}(\prod_{j} Quot(N, r_j, d_j))\to H^{\star}(HQuot(N, {\bf r}, {\bf d}))$$ is surjective. The results of Section $2.1$ imply that in fact the cohomology of $HQuot$ is generated by the Chern classes $c_i(\mathcal A_{j,m})$ for $m\in \{-1,0\}$. Here $\mathcal A_{j,m}=R\pi_{\star} \mathcal F_j(m)$ for $m\geq -1$, and $\pi: \p^1\times HQuot\to HQuot$ is the natural projection. Of course, this change of generators could also be seen more directly. 

When $n\geq 3$, the statement in Lemma $\ref {nomarksflags}$ (bis) follows by the same argument we used in Section $2.2$ for Grassmannians. 

To deal with fewer marked points, we will explain that the map $i^{\star}$ is surjective in $G$-equivariant cohomology. Here $G$ denotes one of the groups $SL_2$, $N$ or $\cs$ which we used in Section $2.3$. Surjectivity is a consequence of the collapse of the Leray spectral sequence in equivariant cohomology as proved for example in \cite {Gi}. Strictly speaking, the group $N$ is not covered by the results of \cite {Gi}, but the argument in the algebraic category presented below takes care of this case as well.

We obtain the diagram:
\begin {center}\xymatrix {
E_{2}^{p,q}=H^{p}_G \otimes H^{q}\left(\prod_{j}Quot(N, r_j, d_j)\right)\ar[d]^{i^{\star}}& \Rightarrow & H^{p+q}_G\left(\prod_{j}Quot(N, r_j, d_j)\right)\ar[d]^{i^{\star}_G}\\
E_{2}^{p,q}=H^{p}_G \otimes H^{q}(HQuot)& \Rightarrow & H^{p+q}_{G}(HQuot)}\end {center}
Surjectivity of the equivariant cohomology restriction map $i_G^{\star}$ follows now from the non-equivariant statement. 
Indeed, it suffices to observe that the collapse of the spectral sequence shows that the equivariant groups on the right admit filtrations such that $i^{\star}$ induces surjections between their associated graded algebras. 

We dedicated the previous subsection to the computation of the equivariant cohomology of the $Quot$ schemes. From the above, we obtain a generation result of the equivariant cohomology of $HQuot$ in terms of the equivariant tautological Chern classes of the bundles $\mathcal A_{j,m}$. Now the arguments which occupy the rest of Section $2.3$ can be used to finish the proof of Lemma $\ref {nomarksflags}$ (bis).

The argument in the Chow groups is slightly more involved, but we will include it here for completeness. For simplicity, let us write $X$ and $Y$ for $HQuot$ and $\prod_{j}Quot(N, r_j, d_j)$ respectively, and then $i:X\to Y$ is the forgetful morphism. We know that the pullback: $$i^{!}: A_{\star}(Y)\cong A^{\star}(Y)\stackrel{i^{\star}}{\longrightarrow}A^{\star}(X)\cong A_{\star}(X) \text{ is surjective}.$$ We want to derive the same statement equivariantly for the action of the groups $SL_2$, $N$ and $\cs$. As before, let $W$ be an open smooth subvariety of an affine space with a free action of $G$. We claim that the map $i^{!}_{W}:A_{\star}(Y\times W)\to A_{\star}(X \times W)$ is also surjective. This follows from the surjectivity of $i^{!}$ and of the following two maps:
$$A_{\star}(Y)\otimes A_{\star}(W)\to A_{\star}(Y \times W)\text { and }  A_{\star}(X) \otimes A_{\star}(W)\to A_{\star}(X \times W).$$ To see that the exterior product maps are surjective, we make use of the fact that both $X$ and $Y$ are smooth projective varieties admitting torus actions with isolated fixed points \cite {quot}. Such torus actions are obtained  from a generic torus action on $\p^1$ and on the fibers of the sheaf $\mathcal O_{\p^1}^{N}$ whose quotients give the $Quot$ schemes. Therefore, the theorem of Bia\l{}ynicki-Birula shows that $X$ and $Y$ can be stratified by unions of affine spaces. For affine spaces surjectivity is clear. Our claim follows inductively, by successively building $X$ and $Y$ from their strata. 

To finish the proof, it suffices to explain the surjectivity of the map: $$i_G^{!}: A^{G}_{\star}(Y)=A_{\star}(Y\times_{G} W)\to A_{\star}(X\times_{G} W)=A^{G}_{\star}(X).$$ We have a fiber diagram: 
\begin {center}
$\begin {CD}
X\times W@>{i_W}>>Y \times W\\
@V{\pi_X}VV @V{\pi_Y} VV\\
X\times_{G} W@>{i_G}>> Y\times_{G}W 
\end {CD}$
\end {center}
where the vertical arrows are principal $G$ bundles. Let $\alpha$ be any class in $A_k(X\times_{G} W)$. Then, our assumption and Theorem $1$ in \cite {V2} respectively show that there are classes $\overline \beta$ and $\beta$ on $Y\times W$ and $Y\times_{G} W$ such that: $$\pi_X^{\star} \alpha = i_W^{!} \overline \beta \; \text { and } \overline \beta=\pi_{Y}^{\star} \beta.$$ Therefore, $$\pi_{X}^{\star} \alpha = i_W^{!} \pi_{Y}^{\star}\beta = \pi_{X}^{\star} i_G^{!} \beta.$$ We first assume $G$ is either $SL_m$ or $GL_m$ and use Theorem $2$ in \cite {V2}. Then, $$\alpha = i_{G}^{!} \beta+\sum_{i\geq 1} c_i^{X} \cap \alpha_i$$ for some classes $\alpha_i\in A_{k+i}(X\times_G W)$ and for some operational classes $c_i^{X}$ on $X\times_{G} W$. Moreover, by \cite {V2}, we find classes $c_i^{Y}$ operating on $A_{\star}(Y \times_G W)$ with $i_G^{\star} c_i^{Y}=c_i^{X}$. Inductively on codimension, we know $\alpha_i=i_G^{!} \beta_i$. The following computation concludes the proof: $$\alpha = i_G^{!} \beta+\sum_{i} i_G^{\star} c_i^{Y} \cap i_G^{!} \beta_i = i_G^{!} \left(\beta+\sum_{i} c_i^{Y} \cap \beta_i\right).$$

There is one remaining case needed for our arguments, namely that of the group $N$. Let $E\to X\times_{GL_2} W$ be the vector bundle associated to the principal $GL_2$ bundle: $$X\times W \to X\times_{GL_2} W.$$ Then $X\times_{SL_2} W$ is total space of the bundle of frames in $\det E$. It is equipped with a projection $\eta: X\times_{SL_2} W \to X\times_{GL_2} W$. More importantly, $X \times_{N} W=(X\times W)\times_{SL_2} SL_2/N$ is the projective bundle $\mathbb P(\eta^{\star} E)$ over $X\times_{SL_2} W$ and comes equipped with a tautological bundle $\tau_{X}$. Then $A_{\star}(X\times_N W)$ is generated by the class $c_1(\tau_X)$ and $A_{\star}(X\times_{SL_2} W)$. The analogous discussion holds for $Y$. From here, the surjectivity of $i_{N}^{!}$ follows from that of $i_{SL_2}^{!}$. 

\subsection {The Cohomology of the Hyper-Quot scheme.} In this section we will prove Lemma $\ref {hq}$ using the diagonal trick of Beauville-Ellingsrud-Stromme \cite {ES}. We will express the class of the diagonal embedding $\Delta: HQuot \hookrightarrow HQuot \times HQuot$ as a combination of classes $\pi_1^\star \alpha \cdot \pi_2^\star \beta$ on $HQuot\times HQuot$, where $\alpha$ and $\beta$ are among the tautological classes listed in Lemma $\ref {hq}$. Here $\pi_1, \pi_2$ are the two projections $HQuot\times HQuot\to HQuot$. This will establish Lemma $\ref {hq}$ completely. Since such arguments are well known, we will only sketch some of the details.

Let $\mathcal K$ be the kernel of the natural sheaf morphism on $\p^1\times HQuot\times HQuot$: $$\bigoplus_{i=1}^{l} \text { Hom }(\pi_1^{\star}\e_i, \pi_2^{\star}\f_i)\to \bigoplus_{i=1}^{l-1} \text {Hom }(\pi_1^{\star}\e_i, \pi_2^{\star}\f_{i+1})\to 0.$$ Let $p:\p^1\times HQuot\times HQuot\to HQuot \times HQuot$ denote the natural projection. It can be shown that $p_{\star} \mathcal K$ is a vector bundle whose rank equals the dimension of $HQuot$, essentially by showing that there are no $H^{1}$'s along the fibers of $p$. In turn, this can be observed via the following argument borrowed from \cite {CF}. Assume we are given two geometric points of $HQuot$: $$0\to E_{\bullet}\to V\otimes \o\to F_{\bullet}\to 0, \;\;0\to E'_{\bullet}\to V\otimes \o \to F'_{\bullet}\to 0.$$ These define a morphism $\p^1\to \p^1\times HQuot\times HQuot$ and we let $K$ be the pullback of $\mathcal K$. We have a natural map from a trivial bundle: $$\text {Hom }(V, V)\otimes \o\to \bigoplus_{i=1}^{l} \text {Hom }(E_i, F'_i)\to \bigoplus_{i=1}^{l-1} \text {Hom }(E_i, F'_{i+1})$$ which factors through $K$. One easily checks that the map $\text {Hom }(V, V)\otimes \o\to K$ is generically surjective. Therefore $H^{1}(\p^1, K)=0$. 

Moreover, a section of $\mathcal K$ is canonically obtained from the natural morphisms: $$\pi_1^{\star} \e_i\to V\otimes \mathcal O\to \pi_2^{\star} \f_i.$$ We also obtain a section of the bundle $p_{\star}\mathcal K$ on $HQuot \times HQuot$. This section vanishes precisely along the diagonal. Therefore $[\Delta]=c_{top}(p_{\star}\mathcal K).$ A Grothendieck Riemann Roch computation expresses the Chern character/classes of $p_{\star} \mathcal K$ as combination of classes $\pi_1^{\star}\alpha\cdot \pi_2^{\star}\beta$ where $\alpha, \beta$ are among the candidates we listed in Lemma $\ref {hq}$, as desired.  

\subsection {Cohomology of fibered products.} To finish the proof of Theorem $\ref {main}$ we need to understand the cohomology of the boundary strata. To this end, we will need to make use of the following result about the cohomology of fibered products.

\addtocounter{lemma}{+1}
\begin {lemma}\label {orbilemma}
Assume there is a fiber square where $p_1$ and $p_2$ are proper morphisms of projective orbifolds with surjective orbifold differentials, and $B$ simply connected:
\begin {center}
$\xymatrix{
{Z=X\times_{B} Y} \ar[r]\ar[d] & Y\ar[d]^{p_2}\\
X\ar[r]^{p_1}& B}
$
\end {center}
The cohomology of $H^{\star}(Z)$ is generated by the image of $H^{\star}(X)\otimes H^{\star}(Y)$.
\end {lemma}

{\bf Proof.} We reformulate the statement as follows. We have a fiber diagram: 
\begin {center}
$\begin {CD}
S'@>{i'}>>S\\
@V{\pi'}VV @V{\pi} VV\\
T'@>{i}>> T 
\end {CD}$
\end {center}
where $T'=B,$ $T=B\times B,$ $S'=Z$, $S=X\times Y$. We observe that $i^{\star}:H^{\star}(T)\to H^{\star}(T')$ is surjective. We want to prove that $(i')^{\star}:H^{\star}(S) \to H^{\star}(S')$ is also surjective. This will follow from a more general argument.  

There are two Leray sequences corresponding to the maps $\pi$ and $\pi'$. Their collapsing is a well known result of Deligne \cite {D}. To apply it we need to know that the differentials of both maps are surjective. Strictly speaking, Deligne works with smooth projective varieties, but his result extends to projective orbifolds (the main ingredient in the proof is the Hard Lefschetz theorem, which holds for orbifolds - see \cite {St}). There are natural morphisms between these spectral sequences: 
\begin {center}
$$\begin {array}{ccc}
H^{p}(T, R^{q}\pi_{\star}\mathbb Q) &\;\Rightarrow &\; H^{p+q}(S)\\
i^{\star}\downarrow & \;&  (i')^{\star}\downarrow\\
H^{p}(T', R^{q}\pi^{'}_{\star}\mathbb Q)&\;\Rightarrow &\;H^{p+q}(S')
\end {array}$$
\end {center}

We claim that the second vertical arrow $(i')^{\star}$ is surjective. We first observe that the first vertical arrow $i^{\star}$ is surjective. Indeed, $i^{\star}:H^{\star}(T)\to H^{\star}(T')$ is surjective. The two local systems given by the direct images $R^{q}\pi_{\star}\mathbb Q$ and $R^{q}{\pi'}_{\star}\mathbb Q$ on $T$, $T'$ are trivial, these spaces being simply connected. Surjectivity of ${i}^{\star}$ follows. Because the two spectral sequences degenerate, there are suitable filtrations $F^{\bullet}$ and $F'^{\bullet}$ of $H^{p+q}(S)$ and $H^{p+q}(S')$, compatible with the map $(i')^{\star}$ such that the map $(i')^{\star}: Gr_{F}^{\bullet}\to Gr_{F'}^{\bullet}$ is surjective. It follows inductively that $(i')^{\star}$ restricted to the successive pieces of the filtrations is also surjective, hence $(i')^{\star}:H^{\star}(S)\to H^{\star}(S')$ is surjective. This completes the proof. 

\subsection {The main result.} Let $X$ be an arbitrary $SL$ flag variety. In this subsection we will conclude the cohomological proof of our main result, Theorem $1$. 

We will apply the lemma proved above to the evaluation maps $ev: \overline M_{0,n+1}(X, \beta)\to X$. Then, the cohomology of the fibered product $\overline M_{0, A\cup \{\bullet\}}(X, \beta_A)\times_{X} \overline M_{0, B\cup \{\star\}} (X, \beta_B)$ is generated by classes coming from each factor. To place ourselves in the context of the lemma, we need to show that the (orbifold) differentials have maximal rank i.e. that these differentials are surjective. Recall from \cite {FP} the construction of $\overline M_{0,n+1}(X,\beta)$. First, one rigidifies the moduli problem. Consider an embedding of $X$ into an $r$ dimensional projective space $\p(V)$, and fix $t=(t_0, \ldots, t_r)$ a basis of $V^{\star}$. We define a moduli space of $t$-rigid stable maps to $\p(V)$. This space parametrizes stable maps $f:C\to \p(V)$ of degree $d$ with $n+1$ markings; in addition to the $n+1$ standard markings of the domain $p_i$, we also fix $d(r+1)$ markings $q_{i,j}$ ($0\leq i\leq r$ and $1\leq j\leq d$) stabilizing the domain curve, such that we have an equality of Cartier divisors $f^{\star} t_i=q_{i,1}+\ldots +q_{i,d}$. 

We consider the closed subscheme of the scheme of $t$-rigid maps to $\p(V)$ which factor through the inclusion of $X$ in $\p(V)$. The moduli space of such $t$-rigid maps $M_t$ is a smooth variety. To get to the Kontsevich moduli space, we need to quotient out the action of a finite group. It suffices to show that the map $ev:M_t \to X$ which evaluates at the last point has maximal rank. This is essentially explained in \cite {FP} and it is quite straightforward. Let $(f, C, p_1,\ldots p_{n+1}, q_{ij})$ be a $t$-rigid map and let $p=f(p_{n+1})$. The differential of the evaluation sends the deformation space of the rigid stable map, $\text{Def}$, to $T_p X$. This map factors as the composition of two surjections $\text {Def}\to {\bf {Def}}\to T_p X$. Here ${\bf {Def}} = H^{0}(f^{\star}TX/T_C(-p_{n+1}))$ is the deformation space of the triple $(f, C, p_{n+1})$. The second map is simply the fiber evaluation. It is surjective because for any genus $0$ stable map $f$, $f^{\star}TX$ is generated by global sections \cite {FP}. 

We are now ready to prove Theorem $\ref {main}$. The statement is proved by double induction, first on the degree $\beta$, and then on the number of markings. In degree $\beta=0$, the result follows from Keel's theorem. Next, we consider the indecomposable classes $\beta$ and $n\leq 1$. When $X=G(k, V)$, the result is a consequence of the description of the moduli space $\overline {\mathcal M}_{0,0}(X, 1)$ as the flag variety $Fl(k-1, k+1, V)$ of two step flags in $V$ of dimensions $k-1$ and $k+1$ respectively (this is explained for example in \cite {STi} Lemma 3.2). For $n=1$, there is a similar description of the moduli space as the flag variety $Fl(k-1, k, k+1, V)$. For general flags $X$, the class $\beta$ is Poincare dual to $c_1(\q_i)$, for some $i$. There exists a flag variety $Y$ (obtained by skipping the $i^{\text {th}}$ quotient in $X$) and a projection morphism $\pi: X\to Y$ such that $\pi_{\star} \beta=0$. All stable maps to $X$ in the class $\beta$ are entirely contained in the fibers of $\pi$ which are Grassmannians. Therefore, the moduli space $\overline {\mathcal M}_{0,0}(X, \beta)$ maps to $Y$, the fibers being flag varieties as above. The main theorem follows immediately.

All other moduli spaces for higher values of $n$ or $\beta$ have nonempty boundary divisors $\overline D(A, B, \beta_A, \beta_B)$, where either $\beta_A$ and $\beta_B$ are both smaller than $\beta$, or $A, B$ have fewer than $n$ points. The cohomology of the fibered product scheme $\overline {M}(A, B, \beta_A, \beta_B)$ (and hence of the fibered product stack) dominating the boundary is computed by Lemma $\ref {orbilemma}$. It is spanned by tautological classes in the light of the induction assumption. We apply the Deligne spectral sequence, specifically the exact sequence $\eqref {exactseq}$. We start with an arbitrary codimension $k$ class $\alpha$ in $\fstack$. Its restriction $j^{\star} \alpha$ on $\mathcal M_{0,n}(X, \beta)$ has Hodge weight $k$. By virtue of Lemma $\ref{nomarksflags}$ (bis), we derive that $j^{\star}\alpha$ is equal to the restriction $j^{\star} \alpha'$ of some tautological class $\alpha'$. Then \eqref{exactseq} shows that $\alpha -\alpha'$ is supported on $\overline {\mathcal M}(A, B, \beta_A, \beta_B)$, hence inductively it is sum of tautological classes. This proves the theorem. 

\section {Some Open Questions} 

The results discussed in this paper should hold for general flag varieties over any field (see item (iii) in the Introduction, under ``Conventions and remarks"). We conjecture:

\begin {conjecture}\label {conj1}
If $X=G/P$ is any generalized variety, where $G$ is a semisimple algebraic group and $P$ a parabolic subgroup, then all rational Chow classes of $\fstack$ are tautological.
\end {conjecture}

We also propose the following conjecture based on computations in \cite {chow} and \cite {O2}:

\begin {conjecture} If $X=G/P$ is any generalized flag variety, (the ranks of) the Chow groups of the open part $\mathcal M_{0,n}(X, \beta)$ are independent on the degree $\beta$, provided that $\beta$ is a linear combination of the indecomposable classes with positive coefficients. 
\end {conjecture}

Since the number of tautological classes is large compared to the rank, it would be interesting to find a universal way of producing relations. For instance, relations can be obtained using Mumford's method: one exploits the vanishing of higher Chern classes of tautological bundles which can otherwise be computed by Grothendieck-Riemann-Roch. 

\begin {question}\label {q1} \begin {itemize}\item [(i)] Find a universal (and complete) way of obtaining relations in the tautological rings. \item [(ii)] Find a concrete description of the cohomology ring when the target is any generalized flag variety.\end {itemize}\end {question} 

In \cite {O2} we show that in codimension $1$, all relations are {\it essentially consequences} of the topological recursion relations coming from ${\overline M}_{0,n}$. This is a manifestation of the fact that the tautological systems are {\it generally} insensitive to the geometry of the target space other than through the evaluation pullbacks. One may hope that all relations between the tautological generators belong to the system of tautological relations preserved by the natural pushforwards and containing, in addition to the usual relations in Gromov-Witten theory, the topological recursion relations. We will make these concepts precise elsewhere.

\begin {question} Is it true that all (non-trivial) relations between the tautological classes (or better, stabilized cohomology classes in the sense of \cite {BH}) are {\it consequences} of the topological recursion relations? 
\end {question}

Finally, one can define the higher genus tautological systems. It is beyond the scope of this paper to discuss this case, but its study in the context of Gromov-Witten theory will be of interest. One possible application of these ideas concerns reconstruction of Gromov-Witten invariants. A genus $0$ example is worked out in our paper \cite {O2} exploiting the relations between the tautological classes.

\appendix

\section {The higher genus tautological systems}

To define the tautological systems $\mathcal R$ for targets $X$ other than $G/P$'s or in higher genera we need to make use of the virtual fundamental classes. As the moduli spaces involved are not necessarily smooth, we will ignore any possible {\it ring} structure $\mathcal R$ may have. 

We let $\mathcal R_{g, n, \beta}\subset A_{\star}(\overline {\mathcal M}_{g,n}(X, \beta))$ be the minimal system satisfying the requirements:
\begin {itemize}
\item $ev_{1}^{\star}\alpha_1\cdot \ldots\cdot ev_{n}^{\star}\alpha_n \cdot \psi_1^{a_1}\cdot \ldots\cdot\psi_{n}^{a_n}\cap \left[\overline {\mathcal M}_{g,n}(X, \beta)\right]^{vir} \in \mathcal R_{g, n, \beta},\text { for } \alpha_i\in A^{\star}(X), \; a_i\geq 0$,
\item the system is closed under pushforward by the forgetful morphisms,
\item the system is closed under the gluing maps.
\end {itemize}

To define the gluing maps $\zeta_{\Gamma}$, we fix a stable modular dual graph $\Gamma$ of genus $g$, degree $\beta$, with $n$ legs. For each vertex $v$, we write $g_v$, $\beta_v$, $n_v$ for the corresponding genus, degree and total valency (the number of incident half edges and legs). The boundary stratum $\mathcal {\overline M}(\Gamma)$ of maps with fixed dual graph $\Gamma$ is obtained from the fibered diagram below, where $E(\Gamma)$ and $H(\Gamma)$ stand for the set of edges and half edges of $\Gamma$, and $\Delta_{\Gamma}$ is the natural diagonal map:

\begin {center}$
\xymatrix{\overline {\mathcal M}_{g,n}(X, \beta)&\overline {\mathcal M}(\Gamma)\ar[l]_{\hspace{.2in}gl_{\Gamma}}\ar[r]\ar[d]&\prod_{v}\overline {\mathcal M}_{g_v, n_v}(X, \beta_v)\ar[d]\\
& X^{E(\Gamma)}\ar[r]^{\Delta_{\Gamma}}& X^{H(\Gamma)}
}$
\end {center}
The gluing map $\zeta_{\Gamma}$ is obtained as composition of the gluing pushforward, the Gysin morphism and the exterior product (which we will omit from the notation): $$\zeta_{\Gamma}:\bigotimes_{v} A_{\star}\left(\overline {\mathcal M}_{g_v, n_v}(X, \beta_v)\right)\to A_{\star}\left(\overline {\mathcal M}_{g,n}(X, \beta)\right),\;\; \zeta_{\Gamma}=(gl_{\Gamma})_{\star} \Delta_{\Gamma}^{!}.$$ By the above definition, we require
$$\zeta_{\Gamma}\left(\bigotimes_{v} \mathcal R_{g_v, n_v, \beta_v}\right)\subset \mathcal R_{g, n, \beta}.$$

\begin {lemmaa} In genus $0$ and for flag varieties, we recover the definition proposed in the Introduction. 
\end {lemmaa} 

The issues are of course closure under multiplication and the $\psi$ classes. To begin, we define the $\kappa$ classes using the forgetful pushforward: $$\kappa_n(\alpha_{n+1}, \ldots, \alpha_{n+p})=\pi_{\star}(ev^{\star}_{n+1}\alpha_{n+1}\cdot \ldots \cdot ev_{n+p}^{\star}\alpha_{n+p}\cap \left[\overline{\mathcal M}_{0,n+p}(X, \beta)\right]^{vir}).$$ 

We let $\mathcal S^d=\mathcal S^{d}_{\beta,n}\subset A_{\star}(\overline {\mathcal M}_{0,n}(X, \beta))$ be the following collection of descendant classes: \begin {equation}\label{desc}\theta=ev_{1}^{\star} \alpha_1 \cdot \ldots \cdot ev_{n}^{\star}\alpha_n\cdot \psi_1^{a_1} \cdot \ldots \cdot \psi_n^{a_n} \cap \kappa_n(\alpha_{n+1}, \ldots, \alpha_{n+p})\in A_{\star}(\overline {\mathcal M}_{0,n}(X, \beta)).\end {equation} We let $\widetilde {\mathcal S}^d$ be the collection of classes: $$\zeta_{\Gamma}(\theta_{\Gamma}), \text { where } \theta_{\Gamma}=\prod_{v} \theta_v, \text { and } \theta_v\in {\mathcal S}^{d}_{\beta_v, n_v},$$ for all stable modular dual graphs $\Gamma$. We define the similar collections of primary classes $\mathcal S^p$ and $\widetilde {\mathcal S}^{p}$ only allowing $a_1=\ldots=a_n=0$ in equation $\eqref{desc}$. 

The lemma will follow from the facts below. 
\begin {itemize}
\item [(a)] For all targets $X$, $\widetilde {\mathcal S}^{p}$ is preserved by the natural pushforwards. 
\item [(b)] For flag vartieties $X$, $\widetilde {\mathcal S}^{d}$ is closed under the multiplication in the Chow ring.  
\item [(c)] For all targets $X$, $\widetilde {\mathcal S}^{d}=\widetilde {\mathcal S}^{p}$; consequently, both collections give additive generators for the tautological systems.
\end {itemize}

To prove $(a)$, we first observe that closure of $\widetilde {\mathcal S}^{p}$ under the gluing pushforwards is obvious. We check closure under the forgetful morphism $\pi:\overline {\mathcal M}_{0,n}(X, \beta)\to \overline {\mathcal M}_{0,n-1}(X, \beta)$. Letting $\Gamma'$ be the graph obtained from $\Gamma$ by forgetting the $n^{\text {th}}$ leg, we obtain: $$\pi_{\star} \zeta_{\Gamma}(\theta_{\Gamma})=\pi_{\star} (gl_{\Gamma})_{\star}\Delta^{!}_{\Gamma} \theta_{\Gamma}=(gl_{\Gamma'})_{\star} \pi_{\star} \Delta_{\Gamma}^{!}\theta_{\Gamma}=(gl_{\Gamma'})_{\star} \Delta_{\Gamma'}^{!} (\pi_{\star}\theta_{\Gamma})=\zeta_{\Gamma'} (\pi_{\star}\theta_{\Gamma}).$$ Special care must be taken when the graph $\Gamma'$ is unstable. At any rate, we reduce our check to classes in $\mathcal S^p$. Then, let $\theta$ be a class as in $\eqref{desc}$, with $a_1=\ldots = a_n=0$. The projection formula shows: $$\pi_{\star}\theta=ev_1^{\star}\alpha_1 \cdot \ldots \cdot ev_{n-1}^{\star}\alpha_{n-1} \cap \kappa_{n-1}(\alpha_n, \ldots, \alpha_{n+p})\in \mathcal S^p.$$

To check $(b)$, we follow an idea of \cite {taut}. We fix two classes $\zeta_{\Gamma} (\theta_{\Gamma})$ and $\zeta_{\Gamma'}(\theta_{\Gamma'})$ supported on $\overline {\mathcal M}(\Gamma)$ and $\overline {\mathcal M}({\Gamma'})$. Their product is computed by the excess intersection formula. The excess bundle will be distributed over the components $\overline {\mathcal M}(\Gamma'')$ of the stack theoretic intersection of $\overline {\mathcal M}(\Gamma)$ and $\overline {\mathcal M}(\Gamma')$. These are indexed by dual graphs which are given additional structure. The graph $\Gamma''$ is endowed with two collapsing maps $\Gamma''\to \Gamma$ and $\Gamma''\to \Gamma'$ which replace whole subgraphs of $\Gamma''$ with vertices of either $\Gamma$ or $\Gamma'$, also collecting the incident legs and the degree labels. Moreover, we require that each half-edge of $\Gamma''$ correspond to a half edge in either $\Gamma$ or in $\Gamma'$. Just as in equation $(11)$ in \cite {taut}, we derive that the excess normal bundle splits as sum of line bundles which are expressed in terms of the cotangent lines. 
The top Chern class of the excess bundle equals $$\prod_{e}(-\psi_v-\psi_w)$$ where $v,w$ are vertices lying on an edge $e$ which ``comes" from both $\Gamma$ and $\Gamma'$. Therefore, the excess intersection formula shows that:  $$\zeta_{\Gamma}(\theta_{\Gamma})\cdot\zeta_{\Gamma'}(\theta_{\Gamma'})=\sum_{\Gamma''}\zeta_{\Gamma''}(\theta_{\Gamma''}).$$ We argue that $\theta_{\Gamma''}$ is an exterior product of classes in $\mathcal S^{d}$. To this end, we observe that:
\begin {itemize}
\item [(i)] Pullback under the gluing morphisms $gl_{\Gamma}$ preserves the evaluation classes $ev^{\star}{\alpha}$ and the $\psi$ classes. 
\item [(ii)] The pullback of a $\kappa$ class under gluing is sum of $\kappa$ classes. 
\item [(iii)] Product of $\kappa$ classes is a $\kappa$ class $$\kappa_n(\alpha_1, \ldots, \alpha_p)\kappa_n(\beta_1, \ldots, \beta_q)=\kappa_n(\alpha_1, \ldots, \alpha_p, \beta_1, \ldots, \beta_q).$$
\end {itemize} 
 
Finally, for the last item on our list, namely $\widetilde {\mathcal S}^{d}\subset \widetilde {\mathcal S}^{p}$, it suffices to show that $\mathcal S^{d}\subset \widetilde {\mathcal S}^{p}$ since $\widetilde {\mathcal S}^p$ is invariant under the gluing morphisms. This will follow if we show that: $$\psi_1\cap\_: \widetilde {\mathcal S}^{p}\to \widetilde {\mathcal S}^{p}.$$ Using the projection formula for the boundary maps, observation $(i)$ and the compatibility of Chern classes with the Gysin morphisms, it suffices to prove that: $$\psi_1\cap\_: \mathcal S^p \to \widetilde {\mathcal S}^{p}.$$ Let $\theta$ be given by $\eqref{desc}$. The projection formula for the morphism $\pi:\overline {\mathcal M}_{0, n+p}(X, \beta)\to \overline {\mathcal M}_{0,n}(X, \beta)$ shows that $\psi_1\cap \theta$ is the $\pi$-pushforward of the class: $$\pi^{\star}\psi_1 \cdot ev_1^{\star} \alpha_1 \cdot \ldots \cdot ev_{n+p}^{\star} \alpha_{n+p} \cap \left[\overline {\mathcal M}_{0, n+p}(X, \beta)\right]^{vir}.$$ 
Using the invariance of $\widetilde {\mathcal S}^{p}$ under the forgetful pushforward by $\pi$ established above, it suffices to prove that this class belongs to $\widetilde {\mathcal S}^{p}$. It is an immediate consequence of the commutation between Gysin morphisms and flat pullbacks that: $$ev_i^{\star}\alpha_i\cap \_: \widetilde {\mathcal S}^p\to \widetilde {\mathcal S}^p.$$ It remains to prove that: $$\pi^{\star}\psi_1 \cap \left[\overline {\mathcal M}_{0, n+p}(X, \beta)\right]^{vir}\in \widetilde {\mathcal S}^{p}.$$ When $X=\pr$, this is a consequence of the equation VI $6.17$ in \cite {M} which expresses $\pi^{\star}\psi_1-\psi_1$ as sum of boundaries, and Lemma 2.2.2 in \cite {divisors} which expresses $\psi_1$ in terms of boundaries, evaluation and operational $\tilde \kappa$ classes. The general case follows pulling back under the closed embedding $i: \overline {\mathcal M}_{0,n+p}(X, \beta)\to \overline {\mathcal M}_{0, n+p}(\pr, i_{\star} \beta)$ induced by a projective embedding of $X$. Then, to finish, we cap with the virtual fundamental class.

\end {document}